\newtheorem{theorem}{Theorem}                   
\newtheorem{lemma}[theorem]{Lemma}              
\newtheorem{corollary}[theorem]{Corollary}      
\newtheorem{claim}[theorem]{Claim}              
\newtheorem{conjecture}[theorem]{Conjecture}    
\newtheorem{fact}[theorem]{Fact}                
\newtheorem{definition}{Definition}    

\newcommand{\indtheorem}{\begin{quote}\begin{theorem}}
\newcommand{\indtheoremend}{\end{theorem}\end{quote}}
\newcommand{\indlemma}{\begin{quote}\begin{lemma}}
\newcommand{\indlemmaend}{\end{lemma}\end{quote}}
\newcommand{\indcor}{\begin{quote}\begin{corollary}}
\newcommand{\indcorend}{\end{corollary}\end{quote}}
\newcommand{\indclaim}{\begin{quote}\begin{claim}}
\newcommand{\indclaimend}{\end{claim}\end{quote}}

\newenvironment{proof}
      {\medskip\noindent{\bf Proof :}\hspace{1em}}
      {\hfill$\Box$\medskip}

\makeatletter
\def\qed{
   {%
      \unskip
      \nobreak \hfil
      \penalty 50               
      \hskip 3em                
      \null \nobreak \hfil
      \qedbox
      \parfillskip=\z@skip
      \finalhyphendemerits=\z@  
      \endgraf                  
   }}
\makeatother
\newcommand{\qedbox}{\ \hfill\rule{1.5ex}{1.5ex}\vspace{\smallskipamount}}

\documentclass[12pt]{article}

\usepackage{amssymb}

\title{Some Remarks on the Jacobian Conjecture and Connections with Hilbert's Irreducibility Theorem}  
\author{
Richard J. Lipton\thanks{Georgia Tech, College of Computing, Atlanta, GA 30332 and Telcordia Research, Morristown, NJ 07960. Email: rjl@cc.gatech.edu. Research supported by NSF grant CCF-0431023}
\and Evangelos Markakis\thanks{Georgia Tech, College of Computing, Atlanta, GA 30332. Email: vangelis@cc.gatech.edu}}      

\date{}

\begin{document}

\maketitle

\begin{abstract}
We have two main results.

\begin{enumerate}
	\item Let $P: \mathbb K^n \rightarrow \mathbb K^n$ be a polynomial map with constant nonzero Jacobian, where $\mathbb K$ is any extension of $\mathbb Q$. Then, $P$ has a polynomial inverse if and only if the range of $P$ contains a cartesian product of $n$ universal Hilbert sets.
	\item Let $P: \mathbb K^2 \rightarrow \mathbb K^2$ be a polynomial map with constant nonzero Jacobian, where $\mathbb K$ is an algebraic number field. Then, $P$ is invertible for ``almost all'' rational integers over $ \mathbb K$. 
\end{enumerate}	
\end{abstract}

\section{Introduction}
                            
The goal of this note is to present some remarks on the famous Jacobian Conjecture. 
Let $P:\mathbb K^n\rightarrow \mathbb K^n$ be a map over a field $\mathbb K$ of characteristic 0, and let $J(P)$ denote the determinant of its Jacobian matrix. We have two main contributions:

\begin{enumerate}
	\item Let $P: \mathbb K^n \rightarrow \mathbb K^n$ be a polynomial map with $J(P)$ identically equal to $1$, where $\mathbb K$ is any finite extension of $\mathbb Q$. Then, $P$ is surjective if and only if $P$ has a polynomial inverse. In fact we prove something stronger: $P$ has a poynomial inverse if and only if the range of $P$ contains a product of universal Hilbert sets, which is much weaker than being onto.
	\item Let $P: \mathbb K^2 \rightarrow \mathbb K^2$ be a polynomial map with $J(P)$ identically equal to $1$, where $\mathbb K$ is an algebraic number field. Then, $P$ is invertible for ``almost all'' integers over $\mathbb K$. 
\end{enumerate}

The first result holds in all dimensions. Two remarks are in order. First, the interesting direction, of course, is the direction from ``sufficiently onto'' to ``invertible'', which is based on the existence of universal Hilbert sets. Second, we recently realized that this result essentially follows from van den Dries and McKenna~\cite{surjective}{Proposition 1.2}. In Section~\ref{statement}, we discuss the similarities and differences between our result and theirs. We would still like to present our proof as it based on a different approach. 

The second result is only proved in the case of two dimensions. There are two main ingredients in our proof. The first is the use of quantitative forms of Hilbert's irreducibility theorem, i.e., counting the number of integers in a certain interval that preserve irreducibility of polynomials. The second tool, which is also the reason that our proof holds only in two dimensions is a characterization of the class of automorphisms that have finite order. In two dimensions it is known that these are all conjugate to linear maps and this is essential to our proof. Extending our techniques beyond two dimensions requires that we understand the structure of automorphisms with finite order 
in higher dimensions, currently an open problem.


As explained above, one of our main tools in both results is the use of various forms of Hilbert's irreducibility theorem and its implications. We believe that the connection between this theorem and the invertibility of polynomial maps is worth further investigation.

\section{Definitions and Basic Facts}

We first state some basic definitions and results that we need.
Suppose that $P: \mathbb K^n \rightarrow \mathbb K^n$ is a polynomial map where $ \mathbb K$ is a field of characteristic zero. This means, as usual, that $P=(f_1,\dots,f_n)$ and each $f_i$ is in $\mathbb K[x_1,\dots,x_n]$. We use $J(P)$ to denote the determinant of the Jacobian matrix $(\frac{\partial F_i}{\partial x_j})_{1\leq i,j\leq n}$ of the map $P$. The famous Jacobian Conjecture states that if $J(P) \equiv 1$, then $P$ has a polynomial inverse.

Let $P=(f_1,\dots,f_n): \mathbb C^n \rightarrow \mathbb C^n$ be a polynomial map with $J(P)\equiv 1$. The following basic facts will be used later on:

\begin{lemma}
\label{independent} 
The functions $f_1,\dots,f_n$ are algebraically independent over $\mathbb C$.
\end{lemma}
\begin{proof}
See~\cite{essen} {Proposition 1.1.31}.
\end{proof}

\begin{lemma}
\label{algebraic}
	Each of $x_1,\dots,x_n$ is algebraic over $\mathbb Q[f_1,\dots,f_n]$.
\end{lemma}

Suppose that $\Phi(u_1,\dots,u_n,z)$ is a polynomial. We say that it {\it depends} on $z$ provided that when written as a polynomial in $z$, i.e., as
	\[ a_m(u_1,\dots,u_n)z^m+\dots+a_0(u_1,\dots,u_n) \]
then $m>0$ and the polynomial $a_m(u_1,\dots,u_n)$ is nonzero.  
	
Lemma~\ref{algebraic} implies the following:
\begin{lemma}
\label{phi}
	For each $x_i$, $i=1,\dots,n$, there is an irreducible polynomial $\Phi_i(u_1,\dots,u_n,z)$ with integer coefficients so that $\Phi_i(f_1,\dots,f_n,x_i)=0$. Moreover, each $\Phi_i(u_1,\dots,u_n,z)$ depends on $z$.
\end{lemma}

The next lemma holds for any dimension, however we will need it only for 
$2$-dimensional maps.

\begin{lemma}
\label{fibers}
	Let $P: \mathbb C^2 \rightarrow \mathbb C^2$ be a polynomial map with $J(P) \equiv 1$. Then, for each $(a,b)\in\mathbb C^2$, $|P^{-1}(a,b)|$ is finite.
\end{lemma}
\begin{proof}
	Let $P=(f_1,f_2)$ such that $J(P) \equiv 1$, and let $a$ and $b$ be given. We need to show that there is only a finite number of solutions to the equations:
\begin{eqnarray}
	 f(x,y)=a \\
	 g(x,y)=b.
\end{eqnarray}
If the polynomials $f(x,y)-a$ and $g(x,y)-b$ have no common factor, then by Bezout's theorem, the number of $x$ and $y$ that satisfy the above equations is bounded by the product of the degrees of $f$ and $g$. Hence suppose that these two polynomials have a common factor. Note that the Jacobian of the map $P' = (f-a, g-b)$ is the same as $J(P)$, thus $J(P') = 1$ for every $z\in\mathbb C^2$. But a simple calculation shows that a pair of polynomials with a nonzero constant Jacobian cannot have a common factor. 
\end{proof}


We use $P \circ Q$ to denote as usual the functional composition of two maps $P$ and $Q$. Thus, for any $z \in \mathbb K^n$, $(P \circ Q)(z) = P(Q(z))$.

\begin{fact}
\label{product}
	Let $P$ and $Q$ be polynomial maps from $\mathbb C^n$ to $\mathbb C^n$. Then, 
	\[J(P \circ Q) = J(P)J(Q).\]
\end{fact}

\section{Statement of Main Results}
\label{statement}

In this Section, we state our main results. Let $\mathbb K$ be a field. We consider polynomial maps $P$ from $\mathbb K^n$ to $\mathbb K^n$ that satisfy tha jacobian condition, i.e., $J(P) \equiv 1$. 
\begin{definition}
An infinite set $H \subseteq \mathbb K$ is called a universal Hilbert
set of order $n$ if for any irreducible polynomial $f(u,x_1,\dots,x_n)$, the set of $a$
for which $f(a,x_1,\dots,x_n)$ is reducible, is a finite subset of $H$.
\end{definition}
Hilbert's irreducibility theorem, see e.g.~\cite{red}, implies that universal Hilbert sets exist for any finite extension $\mathbb K$ of $\mathbb Q$ and they can be quite ``thin''. See~\cite{red} for results on constructing Hilbert sets. 

Our first result shows that if $P$ is "sufficiently onto", then $P$ has a polynomial inverse.
\begin{theorem}
\label{theorem1}
Let $P:\mathbb K^n\rightarrow\mathbb K^n$, where $\mathbb K$ is any finite extension of $\mathbb Q$ and $P$ satisfies $J(P)\equiv 1$. If $P(\mathbb K^n) \supseteq H_1 \times H_2\times\dots \times H_n$, for some universal Hilbert sets $H_1, H_2,\dots,H_n$ of order $n$, then $P$ has a polynomial inverse. 
\end{theorem}

Note that the condition that the range of $P$ only contains $H_1 \times\dots\times H_n$ is much weaker than onto. Note also that our result yields an equivalence between being sufficiently onto and being invertible since the reverse direction of Theorem~\ref{theorem1} is trivial. Our proof works in two steps. We first show that $P$ has a rational inverse. Then, as proved by Keller~\cite{keller}, if $J(P)\equiv 1$ and $P$ has a rational inverse, $P$ in fact has a polynomial inverse.
We recently found out that our first step essentially follows from van den Dries and McKenna~\cite{surjective}{Proposition 1.2} (our condition on the range of $P$ implies that the range is, as in their terminology, Hilbert-dense). Their proof is based on a compactness argument similar in spirit to Gilmore and Robinson~\cite{GR}. We would still like to present our proof as we think it is different and based on more elementary arguments.


In our second main result we use the notion of being invertible for
"almost all" elements of a set. We will say that a set $S\subseteq \mathbb Z^2$ contains 
almost all rational integers of $\mathbb K$ if for all large enough $N$, the complement of $S$ satisfies:
$$|\bar{S} \cap [-N,N]\times [-N,N]| = o(N^2)$$ We can similarly define what it means for a property $\Pi$ to hold for almost all integers.
In particular, we will say that a
map $P$ is injective for almost all integers if $P$ is injective on a set $S$
that contains almost all integers, i.e., for $(x,y)\in S$ and
$(x',y')\in S$, $P(x,y)= P(x',y')$ implies that $(x,y)=
(x',y')$.

\begin{theorem}
\label{theorem2}
Let $P: \mathbb K^2 \rightarrow \mathbb K^2$, where $\mathbb K$ is an algebraic number field and $P$ satisfies $J(P)\equiv 1$. Then $P$ is injective for almost all rational integers of $\mathbb K$. 
\end{theorem}

As usual, we use the term rational integers to distinguish $\mathbb Z$ from the set of algebraic integers over $\mathbb K$. The proof of Theorem~\ref{theorem2} is more involved and uses quantitative versions of Hilbert's irreducibility theorem, i.e., estimates on the number of integers within a certain interval that preserve irreducibility of polynomials. Another essential tool in our proof is a result on the structure of $2$-dimensional automorphisms with finite order (Lemma~\ref{conjugate}). An analogous result in higher dimensions would allow us to prove a more general theorem.

Finally we would like to observe that the starting point in both of our results is the use of Lemma~\ref{phi} and various forms or implications of Hilbert's irreducibility theorem. We believe that the connection between invertibility of polynomial maps and irreducibility questions should be further explored.
\section{Proof of Theorem \ref{theorem1}}

\begin{proof}
We prove the theorem for $\mathbb K$ equal to the rationals and for
$n=2$. The general case is similar. Let $P = (f,g): \mathbb Q^2
\rightarrow \mathbb Q^2$ be a map with $J(P) \equiv 1$. By
Lemma~\ref{phi}, there is an irreducible polynomial $\Phi_1$ such that
$\Phi_1(f, g, x_1) = 0$ (similarly a polynomial $\Phi_2$ for
$x_2$). Let
$$
\Phi_1(f,g,x_1) = a_m(f,g)x_1^m + \dots + a_0(f,g)
$$
Lemma~\ref{phi} implies that $m>0$. We claim that there is a choice of rational values $\alpha\in
H_1,\beta\in H_2$ for $f$ and $g$ (in fact there is an infinite number
of such values), such that the polynomial $\Phi_1'(x_1) \equiv
\Phi_1(\alpha,\beta,x_1)\in\mathbb Q[x_1]$ is irreducible over
$\mathbb Q$, it has a rational root and it has the same degree in
$x_1$ as the original $\Phi_1$. To see this, note that for any pair
$(\alpha,\beta) = (f(x_1,x_2),g(x_1,x_2))$, for $(x_1,x_2)\in\mathbb
Q^2$, it is true that $x_1$ is a rational root of
$\Phi_1(\alpha,\beta,x_1)$ and $x_2$ is a root of
$\Phi_2(\alpha,\beta,x_2)$. Suppose we first substitute $f$ with
$\alpha\in H_1$ in $\Phi_1$. By the definition of a Hilbert set, there
is only a finite number of $\alpha$'s that make $\Phi_1(\alpha,g,x_1)$
reducible. Furthermore, there is only a finite number of $\alpha$'s
that make $a_m(\alpha,g)$ identically $0$. Once we fix $\alpha$, then
again there can be at most a finite number of choices for $\beta$ that
either make $\Phi_1(\alpha,\beta,x_1)$ reducible or make the highest
degree term in $x_1$ vanish. Since the range of $P$ contains
$H_1\times H_2$, we can always find a pair $(\alpha,\beta)$ with the
desirable properties. However, if $deg_{\Phi_1}(x_1)>1$, then we have
a contradiction, since $\Phi_1'(x_1)$ is irreducible over $\mathbb Q$
and we have assumed that it has a rational root. The same is true if
$deg_{\Phi_2}(x_2)>1$. Hence $deg_{\Phi_1}(x_1) = deg_{\Phi_2}(x_2) =
1$. Then $x_1, x_2\in\mathbb Q(f,g)$, which means that $P$ has a
rational inverse. Since $J(P)\equiv 1$, it follows by~\cite{keller}
that $P$ in fact has a polynomial inverse.
\end{proof}

\section{Proof of Theorem \ref{theorem2}}
Let $P = (f,g)$. For ease of notation, we use $x$ and $y$ instead of $x_1$ and $x_2$ for the variables on which $f$ and $g$ depend on. We present the proof with $\mathbb K = \mathbb Q$. The generalization to any number field is straightforward.
By Lemma~\ref{phi}, there is a polynomial $\Phi(u_1,u_2,z)$ that depends on $z$, such that $\Phi(f,g,x) = 0$. Similarly there is a polynomial $\Psi(u_1,u_2,z)$ for which $\Psi(f,g,y) = 0$.

Suppose that in $\Phi(u_1,u_2,z)$ we substitute $u_1$ and $u_2$ by $f$ and $g$. We can then see $\Phi$ as a polynomial in $z$ with coefficients from $\mathbb Q[f,g]$:
$$ \Phi(f,g,z) = a_m(f,g)z^m + a_{m-1}(f,g)z^{m-1}+...+a_0(f,g)$$

We can further substitute $f$ and $g$ as functions of $x$ and $y$ and factor the resulting polynomial over $\mathbb Q(x,y)$. We will then obtain a polynomial in 
$\mathbb Q(x,y)[z]$:
\begin{equation}
\label{poly:fgz}
\Phi(f,g,z) = (z - \phi_1(x,y))(z - \phi_2(x,y))...(z - \phi_k(x,y))A(x,y,z)
\end{equation}
where the $\phi_i$'s are rational functions of $x$ and $y$ and $A$ is an irreducible polynomial. We can also assume that each $\phi_i$ has integer coefficients.

Similarly for the polynomial $\Psi$ we have:
$$\Psi(f,g,z) = (z - \psi_1(x,y))(z - \psi_2(x,y))...(z - \psi_l(x,y))B(x,y,z)$$
Note that both polynomials have at least one factor, i.e., $k, l\geq 1$ because $x$ (resp. y) is a root (since $\Phi(f,g,x) = 0$). 

Let $u,v$ be the values of $f$ and $g$ at a point, say $u = f(x_0,y_0)$ and $v = g(x_0,y_0)$, for some $(x_0,y_0)\in\mathbb Q^2$. We want to see when can we say that the pair $(u,v)$ has no other preimage. We will show that there exists a set $S$ that contains almost all integer pairs, such that for any $(x_0,y_0)\in S$, the corresponding pair of values $(u,v)$ has no other preimage within that set.

From now on, we assume that $(x_0,y_0)\in \mathbb Z^2$ and that both $x_0$ and $y_0$ are in $[-N,N]$, for some large enough $N$. Throughout our proof, we will eliminate integer pairs from $[-N,N]^2$ for which our arguments do not apply. We call such pairs "bad" pairs. We will show that there is a constant $n_0$ such that for all $N\geq n_0$, the number of bad pairs is $o(N^2)$. This will directly imply that the map $P$ is injective on a set that contains almost all integer pairs.  

Substituting $(x_0,y_0)$
in $\Phi, \Psi$ would yield the following two univariate polynomials:

\begin{eqnarray}
\Phi(z) & = & (z-\alpha_1)...(z-\alpha_k)A(z)\\
\Psi(z) & = & (z-\beta_1)...(z-\beta_l)B(z)
\end{eqnarray}
where $\alpha_i = \phi_i(x_0,y_0), ~i =1,\dots,k$, $\beta_j = \psi_j(x_0,y_0),~j=1,\dots$, $A(z) = A(x_0,y_0,z)$ and $B(z) = B(x_0,y_0,z)$.

We first note that for almost all integer pairs $(x_0,y_0)$, the polynomials $A(z), B(z)$ are
irreducible over $\mathbb Q$, which follows from the result of~\cite{cohen}, a quantitative form of Hilbert's irreducibility theorem.
In particular, if we
substitute $x,y$ with integer values in the interval $[-N, N]$, there
can be at most $O(N^{3/2}\log{N})$ bad pairs that make $A(x,y,z)$ reducible out of a
total of $O(N^2)$ possible pairs (see~\cite{red}{$~$Chapter 4} for related results). 

Consider an integer pair $(x_0,y_0)$ such that $A(x_0,y_0,z)$ and $B(x_0,y_0,z)$ are irreducible over $\mathbb Q$. Then the only rational roots of $\Phi(z), \Psi(z)$ are the $\alpha$'s and
the $\beta$'s. Notice also that for all the preimages of $(u,v)$, say
$\{(x_i,y_i)\}$, it holds that $x_i$ is a root of $\Phi(z)$ and $y_i$
is a root of $\Psi(z)$. This comes from the fact that $\Phi$
and $\Psi$ satisfy $\Phi(f,g,x) = 0$ and $\Psi(f,g,y) = 0$. Hence, there is at least one pair, say
$(\alpha_1,\beta_1)$, that
is equal to $(x_0,y_0)$. To see if $u,v$ has any other integer preimage, we only need to do the
following: For every integer root $\alpha_i$ and every integer root $\beta_j$, we check to see whether $P(\alpha_i,\beta_j) = (u,v)$. 
If we find only one such pair, then
$(u,v)$ has no other integer preimage. 

Suppose that for at least two distinct pairs say $(\alpha_1,\beta_1), (\alpha_{i},\beta_{j})$, we get the value $(u,v)$. We claim that this cannot happen for a lot of
integers. One of the two pairs, namely $(\alpha_i,\beta_j)$, is not $(x_0,y_0)$ and is equal to $(\phi_i(x_0,y_0)$ , $\psi_j(x_0,y_0)$. We also have that $P(\phi_i(x_0,y_0)$ , $\psi_j(x_0,y_0)) = P(x_0,y_0)$. Let $Q$ be the bivariate map $Q = (\phi_i,\psi_j)$. Obviously $Q$ is not the identity map. We first show that for almost all integers, we may assume that the map $Q$ is in fact a polynomial map. For this we use the following lemma:
\begin{lemma}
\label{rational}
Let $a(x,y)/b(x,y)$ be a rational function with integer coefficients. Then for large enough $N$, the number of integer pairs $(x_0,y_0)\in [-N,N]^2$ for which $a(x_0,y_0)/b(x_0,y_0)$ is an integer is $O(N^{3/2})$.
\end{lemma}
\begin{proof}
Assume without loss of generality that $b$ is irreducible. We estimate separately for each $y_0\in [-N,N]$, the number of $x_0$'s such that
$b(x_0,y_0)$ divides $a(x_0,y_0)$. There are two cases to consider for $y_0$. First suppose that $b(x,y_0)$ becomes reducible. The result of Fried~\cite{fried}, which is a consequence of Hilbert's irreducibility theorem, implies that there can be at most $O(\sqrt{N})$ such $y_0$'s. Hence there can be at most $O(N^{3/2})$ such pairs $(x_0,y_0)$ for which the rational function takes an integer value. Assume now that $b(x,y_0)$ remains irreducible, which happens for $O(N)$ values of $y_0$. Let $R(y_0)$ be the resultant of $a(x,y_0)$ and $b(x,y_0)$, which is a polynomial in $y_0$ and let $r$ be the degree of $R(y_0)$ (for a definition of the resultant, see ~\cite{lang}). We consider two subcases. Suppose that $R(y_0) = 0$. This can happen for at most $r$ values of $y_0$ and by picking $N$ large enough, we can make $r$ as small as $N^{\epsilon}$ for any small $\epsilon >0$. Therefore there can be at most $O(N^{1+\epsilon})$ pairs $(x_0,y_0)$ that fall under this subcase. Assume now that $R(y_0)\neq 0$, which is true for $O(N)$ values of $y_0$. This implies that $a(x,y_0)$ and $b(x,y_0)$ are relatively prime and there is a $d\in\mathbb Z$ and polynomials $q, s\in\mathbb Z[x]$ such that: $$q(x)a(x,y_0)+ s(x)b(x,y_0) = d$$

For $b(x_0,y_0)$ to divide $a(x_0,y_0)$, it has to be the case that $b(x_0,y_0)$ is equal to a divisor of $d$ (or minus a divisor of $d$). However for any $\delta >0$, the number of divisors of any large enough number $n$ is $O(n^{\delta})$~\cite{HW}. We also know that $d$ is at most a polynomial in $y_0$ (assume that $y_0>0$) by the way it was constructed and therefore for any $\epsilon>0$, we can choose large enough $N$ so that there are at most $O(N^{\epsilon /2})$ divisors of $d$. Then for each $(x_0,y_0)$ that we are interested in, $x_0$ has to be a solution to $b(x,y_0) = d'$ for some divisor $d'$ of $d$. Hence, for each $y_0$ in this subcase, we can have at most $O(deg(b)N^{\epsilon /2})$ values for $x_0$ that make $b(x_0,y_0)$ divide $a(x_0,y_0)$. Therefore the total number of pairs $(x_0,y_0)$ can be made $O(N^{1+\epsilon})$. Finally, summing up all the integer pairs that we counted in each case, we get a total of $O(N^{3/2})$. 
\end{proof}

Therefore, even if all the $\phi_i$'s and $\psi_j$'s were rational functions, there can be at most $O(klN^{3/2})$ integer pairs that achieve integer values under some $\phi_i$ and $\psi_j$. By choosing $N$ large enough we can ensure that this is only $o(N^2)$. Since in our test we discard noninteger $\alpha$'s and $\beta$'s, then for almost all integer pairs $(x_0,y_0)\in [-N,N]^2$, the values $\alpha_i$ and $\beta_j$ as described above, are both integers only if the functions $\phi_i$ and $\psi_j$ are polynomials. Hence, $Q$ is a polynomial map satisfying $(P\circ Q)(x_0,y_0) = P(x_0,y_0)$. We consider the following two cases:\\
         
\noindent {\bf Case 1} $P\circ Q$ is not identical to $P$. In this case, the equation $P(Q(x,y)) - P(x,y) = 0$ defines a non-trivial variety. But varieties can hit only a small fraction of integers in $[-N,N]^2$ as guaranteed by the following well known lemma (e.g., see~\cite{red}{$~$Lemma 1, p. 298):
\begin{lemma}
Let $V$ be the variety: $V = \{(x,y): R(x,y) = 0, R\in\mathbb Q[x,y]\}$. For any $\epsilon >0$, there is large enough $N$ so that the number of integer points in $[-N,N]^2$ that belong to $V$ is $O(N^{1+\epsilon})$.
\end{lemma} 
 
\noindent {\bf Case 2} $P\circ Q\equiv P$
This case is more complicated. We will derive a contradiction by showing that $Q$ has to be the identity map.
First note that since $P(Q)\equiv P$ over $\mathbb Q$, the same will hold over $\mathbb C$. From now on we look at $P$ and $Q$ as polynomial maps from $\mathbb C^2$ to $\mathbb C^2$. Note also that by Fact~\ref{product}, we have $J(Q) \equiv 1$. 

We use $P^t$ to denote the $t$-fold composition of the map $P$ with itself. Thus, $P^2 = P \circ P$. We say that $P$ is {\it conjugate} to a linear map if there exists an invertible polynomial map $S$ and a linear map $L$ so that $P = S^{-1} \circ L \circ S$. As usual a map $P=(f,g)$  is {\it linear} provided each polynomial $f$ and $g$ has degree at most one.

In the rest of our analysis, we make repeated use of the following lemma:
\begin{lemma}
\label{fixedpoint}
If the map $Q$ has a fixed point, then $Q$ is the identity map.
\end{lemma}
\begin{proof}
The proof is based on the inverse function theorem. Suppose $Q$ has a fixed point, say $Q(a,b) = (a,b)$, where $(a,b)\in\mathbb C^2$. By the inverse function theorem, we know that $P$ is locally invertible at a neighborhood of $(a,b)$, i.e., there exists an open set $U$ containing $(a,b)$ and an open set $V$ containing $P(a,b)$, such that $V = P(U)$ and $P$ is one-to-one, when restricted to $U$.  We can pick a small enough open subset of $U$, say $D\subseteq U$, such that for every $(a',b')\in D$, $Q(a',b')\in U$. Since $P(Q)\equiv P$, we have that 
$$ P(Q(a',b')) = P(a',b') ~~\forall (a',b')\in D$$
But $P$ is one-to-one, when restricted to $U$. It follows that $Q(a',b') = (a',b')$ on the open set $D$ and since $Q$ is a polynomial map, this implies that $Q$ has to be the identity map.
\end{proof}

Lemma~\ref{fixedpoint} enables us to prove the following property of the map $Q$.
\begin{lemma}
\label{order}
The map $Q$ has a finite order, i.e., there exists a positive integer $t\geq 2$ such that $Q^t$ is the identity map.
\end{lemma}
\begin{proof}
Pick $z\in\mathbb C^2$ and let $u = P(z)$. Consider the terms $z, Q(z), Q^2(z)$,$\dots$. By Lemma~\ref{fibers} we know that $|P^{-1}(z)|$ is finite. On the other hand, $u = P(z) = P(Q(z)) = P(Q^2(z)) = \dots$. Hence there exist $r>s$ such that $Q^r(z) = Q^s(z)$. This means that the map $Q^{r-s}$ has a fixed point and it also satisfies $J(Q^{r-s}) \equiv 1$ and $P\circ Q^{r-s}  = P$. Lemma~\ref{fixedpoint} completes the proof. Note that since we have assumed that $Q$ is not the identity map, then $r-s\geq 2$. 
\end{proof}

The following lemma, proved in~\cite{kambayashi} (see also~\cite{essen}), is essential in our proof and is the only step in which we need to assume that $P$ is $2$-dimensional. An analogous result in higher dimensions, currently not known to the best of our knowledge, would imply a generalization of our result as well.
\begin{lemma}{\cite{kambayashi}}
\label{conjugate}
Let $P: \mathbb C^2 \rightarrow \mathbb C^2$ be a polynomial map for which there is a $t\in\mathbb N$ with $P^t = I$. Then $P$ is conjugate to a linear map.
\end{lemma}

\begin{lemma}
\label{fixed}
Let $Q: \mathbb C^2 \rightarrow \mathbb C^2$ be a polynomial map with $J(Q) \equiv 1$. Suppose further that $Q$ has finite order. Then, $Q$ has a fixed point.
\end{lemma}
\begin{proof}
Consider a map $Q$ that has finite order, satisfies $J(Q) \equiv 1$ and has no fixed points. We will obtain a contradiction. First note that $Q$ is an invertible polynomial map in two dimensions (tame automorphism) and of finite order. By Lemma~\ref{conjugate}, $Q$ is {\em linearizable}, i.e., there exists an automorphism $R$ such that $R^{-1}QR$ is a linear map. Without loss of generality we can therefore assume that $Q$ itself is a linear map. This is because one can easily check that if $Q$ satisfies all the above properties (no fixed points, finite order and $J(Q) \equiv 1$), then $R^{-1}QR$ also satisfies them. Hence $Q$ is of the form $Q(x,y) = L(x,y)^T+b$, where $L$ is a $2\times 2$ matrix and $b$ is a $2$-dimensional vector. We know that $L$ can be written as $S^{-1}TS$, where $T$ is in Jordan form. Consider the map $SQS^{-1}$. It can be checked again that if $Q$ has all the assumed properties, so does $SQS^{-1}$. But 
$$SQS^{-1}(x,y) = SLS^{-1}(x,y) + Sb = T(x,y)+Sb$$ Therefore, we can assume that the map $Q$ is of the form $T(x,y) + b$, where $T$ is in Jordan form. This means that $T(2,1) = 0$, $T(1,2) = \mu$, where $\mu$ is either $0$ or $1$. Let $T(1,1) = \lambda_1$ and $T(2,2) = \lambda_2$. Notice that since $J(Q) \equiv 1$, we have that $\lambda_1 \lambda_2 = 1$. 

There are 2 cases to consider. First suppose that $\mu = 0$. We know that either $\lambda_1 = \lambda_2 = 1$ or $\lambda_1\neq\lambda_2\neq 1$. The first case implies either that $Q$ is the identity map or that $Q$ does not have a finite order. To see this, note that for any $t\in \mathbb N$, $Q^t(x,y) = (x+tb_1,y + tb_2)$. If the vector $b$ is the $0$-vector, then $Q$ is the identity map and has a fixed point. Otherwise $Q$ cannot have a finite order, a contradiction. In the case that $\lambda_1\neq\lambda_2$ it is easy to show that $Q$ has a fixed point.

In the case that $\mu = 1$, we can do a similar analysis and obtain that $Q$ has a fixed point. 
\end{proof}

Hence $Q$ has a fixed point, which by Lemma~\ref{fixedpoint} implies that $Q$ is the identity map, a contradiction. Therefore $(x_0,y_0)$ is the only preimage of $(u,v)$. In various steps of our analysis we only needed to ignore integer pairs in $[-N,N]^2$ that were no more than $o(N^2)$. Hence the proof of Theorem~\ref{theorem2} is complete.

\section{Conclusions}
We have obtained some connections between Hilbert's irreducibility theorem (in various forms) and invertibility of polynomial maps over algebraic number fields.
Another essential ingredient of our proof is the linearization of polynomial automorphisms with finite order. 

We think it is possible to generalize Theorem~\ref{theorem2} and show that $P$ is injective for almost all algebraic integers over $\mathbb K$. One of the steps that requires a different analysis towards this is Lemma~\ref{rational}.
Another way to enlarge the set on which $P$ is injective in the statement of Theorem~\ref{theorem2} could be to start with a complete factorization of the polynomial $\Phi(f,g,z)$, in which the functions $\phi_i$ and $\psi_j$ would be algebraic functions of $x$ and $y$ and perform a similar analysis. In fact we believe that the jacobian conjecture is equivalent to a statement analogous to Lemma~\ref{fixed}:

\begin{conjecture}
The jacobian conjecture is equivalent to proving the following statement:
Let $P:\mathbb C^2\rightarrow \mathbb C^2$ with:
\begin{enumerate}
\item $J(P) \equiv 1$,
\item there exists an algebraic function defined on some open set ${\cal U}$, such that $P\circ Q = P$,
\item $Q$ has finite order.
\end{enumerate}
Then the map $Q$ has a fixed point.
\end{conjecture}


\end{document}